\documentclass[titlepage,twoside,12pt]{article}
\usepackage{amssymb}
\usepackage{amsfonts}
\begin{document}

{\LARGE \bf Two Generalizations of Tensor \\ \\ Products, Beyond Vector Spaces} \\ \\

% arXiv:physics/0701116v4 corrected on 3/7/8

{\bf Elem\'{e}r E ~Rosinger} \\ \\
{\small \it Department of Mathematics \\ and Applied Mathematics} \\
{\small \it University of Pretoria} \\
{\small \it Pretoria} \\
{\small \it 0002 South Africa} \\
{\small \it eerosinger@hotmail.com} \\ \\

{\bf Abstract} \\

Two successive generalizations of the usual tensor products are given. One can be constructed for arbitrary binary
operations, and not only for semigroups, groups or vector spaces. The second one, still more general, is constructed for
arbitrary {\it generators} on sets. \\ \\

{\large \bf PART 1 : Tensor Products Beyond Vector Spaces} \\ \\

{\bf 1. Preliminaries} \\

Tensor products have lately achieved a significantly greater importance in view of their role in Quantum Mechanics, and
specifically, in Quantum Computation, where entanglement proves to be one of the fundamental assets that make such
computations far faster than on the usual electronic digital computers. In view of that, it is of interest to look deeper
into the mathematical structures which underlie the usual tensor products of vector spaces, and find the simplest and most
general ones which are indeed essential. Here, two successive generalizations of the usual tensor product of vector spaces
are presented. The first one shows that tensor products $X \bigotimes Y$ can in fact be defined for arbitrary sets $X$ and
$Y$, and arbitrary binary operations $\alpha$ and $\beta$ on $X$, respectively, $Y$. This alone is a significantly large
extension, since in the usual case the respective binary operations are restricted to the usual addition of vectors in
vector spaces. The second generalization is still more considerable, since instead of binary operations, one can use the
far more general concept of {\it generators}. \\ \\

{\bf 2. Tensor Products beyond Vector Spaces : \\
        \hspace*{0.45cm} the case of Binary Operations} \\

Let us present the {\it first extension} of the standard definition of {\it tensor product} $X \bigotimes Y$, see Appendix,
to the case of  two structures $( X, \alpha )$ and $( Y, \beta )$, where $\alpha : X \times X \longrightarrow X,~ \beta :
Y \times Y \longrightarrow Y$ are arbitrary binary operations on two arbitrary given sets $X$ and $Y$, respectively. The
way to proceed is as follows. Let us denote by $Z$ the set of all finite sequences of pairs \\

(2.1)~~~ $ ( x_1, y_1 ), \dots , ( x_n, y_n ) $ \\

where $n \geq 1$, while $x_i \in X,~ y_i \in Y$, with $1 \leq i \leq n$. We define on $Z$ the binary operation $\gamma$
simply by the concatenation of the sequences (2.1). It follows that $\gamma$ is associative, therefore, each  sequence
(2.1) can be written as \\

(2.2)~~~ $ ( x_1, y_1 ), \dots , ( x_n, y_n ) =
                 ( x_1, y_1 ) ~\gamma~ ( x_2, y_2 ) ~\gamma~ \ldots ~\gamma~ ( x_n, y_n ) $ \\

where for $n = 1$, the right hand term is understood to be simply $( x_1, y_1 )$. Obviously, if $X$ or $Y$ have at least two elements, then $\gamma$ is not commutative. \\

Thus we have \\

(2.3)~~~ $ Z = \left \{ ( x_1, y_1 ) ~\gamma~ ( x_2, y_2 ) ~\gamma~ \ldots ~\gamma~ ( x_n, y_n ) ~~
                             \begin{array}{|l}
                              ~ n \geq 1 \\ \\
                              ~ x_i \in X,~ y_i \in Y,~ 1 \leq i \leq n
                             \end{array} \right \} $ \\ \\

which clearly gives \\

(2.4)~~~ $ X \times Y \subseteq Z $ \\

Now we define on $Z$ an equivalence relation $\approx_{\alpha, \beta}$ as follows. Two sequences in (2.1) are equivalent, if and only if they are identical, or each can be obtained from the other by a finite number of applications of the following operations \\

(2.5)~~~ permute pairs $( x_i, y_i )$ within the sequence \\

(2.6)~~~ replace $( \alpha ( x_1, x\,'_1 ) , y_1 ) ~\gamma~ ( x_2, y_2 ) ~\gamma~ \ldots ~\gamma~ ( x_n, y_n )$ \\
         \hspace*{1.3cm} with $( x_1, y_1 ) ~\gamma~ ( x\,'_1, y_1 ) ~\gamma~ ( x_2, y_2 ) ~\gamma~ \ldots ~\gamma~
         ( x_n, y_n )$, or vice-versa \\

(2.7)~~~ replace $( x_1, \beta ( y_1, y\,'_1 ) ) ~\gamma~ ( x_2, y_2 ) ~\gamma~ \ldots ~\gamma~ ( x_n, y_n )$ \\
         \hspace*{1.3cm} with $( x_1, y_1 ) ~\gamma~ ( x_1, y\,'_1 ) ~\gamma~ ( x_2, y_2 ) ~\gamma~ \ldots ~\gamma~
         ( x_n, y_n )$, or vice-versa \\

Let us note that, in view of the rather general related result in Lemma 2.1. at the end of this section, the binary
relation $\approx_{\alpha, \beta}$ defined above on $Z$ is indeed an equivalence relation. \\

Finally, the {\it tensor product} of $( X, \alpha )$ and $( Y, \beta )$ is defined to be the
quotient space \\

(2.8)~~~ $ X \bigotimes_{\alpha, \beta} Y = Z / \approx_{\alpha, \beta} $ \\

with the canonical quotient embedding, see (2.4) \\

(2.9)~~~ $ X \times Y \ni ( x, y ) \longmapsto x \bigotimes_{\alpha, \beta} y \in X \bigotimes_{\alpha, \beta} Y $ \\

where as in the usual case of tensor products, we denote by $x \bigotimes_{\alpha, \beta} y$, or simply  $x \bigotimes y$,
the equivalence class of $( x, y ) \in X \times Y \subseteq Z$. Let us show that the mapping (2.9) is indeed injective. Let therefore $( x, y ),
( x\,', y\,' ) \in X \times Y$ be such that $x \bigotimes_{\alpha, \beta} y = x\,' \bigotimes_{\alpha, \beta} y\,'$. Then obviously $( x, y )
\approx_{\alpha, \beta} ( x\,', y\,' )$. However, in view of (2.5) - (2.7), this can only happen when $( x, y ) = ( x\,', y\,' )$. \\

Obviously, the binary operation $\gamma$ on $Z$ will canonically lead by this quotient operation to a {\it commutative}
and {\it associative} binary operation on $X \bigotimes_{\alpha, \beta} Y$, which for convenience is denoted by the same
$\gamma$, although in view of (2.8), this time it depends on $\alpha$ and $\beta$, thus it should rigorously be written
$\gamma_{\alpha, \beta}$. \\

In this way, the elements of $X \bigotimes_{\alpha, \beta} Y$ are all the expressions \\

(2.10)~~~ $ x_1 \bigotimes_{\alpha, \beta} y_1 ~\gamma~ x_2 \bigotimes_{\alpha, \beta} y_2 ~\gamma~
                                   \ldots ~\gamma~ x_n \bigotimes_{\alpha, \beta} y_n $ \\

with $n \geq 1$ and $x_i \in X,~ y_i \in Y$, for $1 \leq i \leq n$. \\

The customary particular situation is when $X$ and $Y$ are commutative semigroups, groups, or even vector spaces over some field $\mathbb{K}$. In this case $\alpha, \beta$ and $\gamma$ are as usual denoted by +, that is, the sign of addition. \\

It is easy to note that in the construction of tensor products above, it is {\it not} necessary for $( X, \alpha )$ or
$( Y, \beta )$ to be semigroups, let alone groups, or for that matter, vector spaces. Indeed, it is sufficient that
$\alpha$ and $\beta$ are arbitrary binary operations on $X$ and $Y$, respectively, while $X$ and $Y$ can be arbitrary
sets. \\

Also, as seen above, $\alpha$ or $\beta$ need {\it not} be commutative either. However, the resulting tensor product
$X \bigotimes_{\alpha, \beta} Y$, with the respective binary operation $\gamma$, will nevertheless be commutative
and associative. \\

An important fact related to the tensor products defined above is that they have a {\it universality} property which is a
natural generalization of the similar one for usual tensor products, see (A1.6.3). This generalized universality property
is presented in section 4. \\

Above, in showing that $\approx_{\alpha, \beta}$ is an equivalence relation on $Z$, we used the following easy to prove \\

{\bf Lemma 2.1.} \\

Let on a nonvoid set $E$ be given a family $( \equiv_i )_{i \in I}$ of {\it symmetric} binary relations. Further, let us
define on $E$ the binary relation $\approx$ as follows. For $a, b \in E$, we have $a \approx b$, if and only if $a = b$,
or there exists a finite sequence \\

$~~~~~~ a = c_0 \equiv_{i_0} c_1 \equiv_{i_1} c_2 \equiv_{i_2} \ldots \equiv_{i_{n-2}} c_{n-1} \equiv_{i_{n-1}} c_n = y $ \\

where $c_1, \ldots , c_{n-1} \in E$. \\

Then $\approx$ is an {\it equivalence} relation on $E$. \\

{\bf Note 2.1.} \\

For clarity about Lemma 2.1. and its role above, let us recall the following standard facts about equivalence relations on
arbitrary sets. Given a nonvoid set $E$, let ${\cal EQ}$ be the set of all equivalence relations on $E$. As is known, each
such equivalence relation $\equiv$ in ${\cal EQ}$  can be identified with a subset $S_\equiv \subseteq E \times E$ in the
following way : given $a, b \in E$, then $a \equiv b \Longleftrightarrow ( a, b ) \in S_\equiv$. \\
Now, we define a {\it partial order} $\prec$ on the set of equivalence relations ${\cal EQ}$ as follows : for two
equivalence relations $\equiv, ~\equiv\,'$ in ${\cal EQ}$, we have $\equiv ~\prec~ \equiv\,' \Longleftrightarrow S_\equiv
\subseteq S_{\equiv\,'}$, where $S_\equiv, S_{\equiv\,'} \subseteq E \times E$ correspond as defined above, to $\equiv,~
\equiv\,'$, respectively. \\
It is easy to see that, in the sense of this partial $\prec$ order, the usual equality $=$ on $E$ is the {\it smallest}
equivalence relation on $E$, and $S_=$ is in fact $\{ ( a, a ) ~|~ a \in E \}$, which is the {\it diagonal} in $E \times
E$. On the other hand, the {\it largest} equivalence relation on $E$, denoted say by $\equiv_{total}$, and which trivially
makes every $a, b \in E$ equivalent with one another, corresponds to the subset $S_{\equiv_{total}} = E \times E \subseteq
E \times E$. \\

A convenient consequence of the above is that, given any subset $R \subseteq E \times E$, there is always a {\it smallest}
equivalence relation $\equiv_R$ in ${\cal EQ}$, such that $R \subseteq S_{\equiv_R}$. Indeed, let $R^*$ be the
intersection of all sets $S \subseteq E \times E$ which correspond to equivalence relations on $E$, and for which $R
\subseteq S$. Obviously, such sets $S$ always exist, since for instance, we can take $S = S_{\equiv_{total}} = E \times
E$. Let now $\equiv_R$ be the binary relation on $E$ defined, for $a, b \in E$, by $a \equiv_R b \Longleftrightarrow ( a,
b ) \in R^*$. Then it is easy to see that $\equiv_R$ is indeed an equivalence relation on $E$. \\

In the particular case when $R$ is symmetric, that is, for $a, b \in E$, we have \\

$~~~~~~ ( a, b ) \in R ~~~\Longrightarrow~~~ ( b, a ) \in R $ \\

then $\equiv_R$ constructed above is called the {\it transitive closure} of R. \\

Now it is easy to note that Lemma 2.1. above does in fact construct such a transitive closure of the family
$( \equiv_i )_{i \in I}$ of symmetric binary relations. \\

Finally, the fact that $\approx_{\alpha, \beta}$ constructed above is indeed an equivalence relation on $Z$, is an
application of transitive closure to the symmetric binary relations on $Z$ defined by (2.5) - (2.7). \\ \\

{\bf 3. Tensor Products beyond Vector Spaces : \\
        \hspace*{0.45cm} the case of Generators} \\

The {\it second extension} of the standard definition of {\it tensor products}, see Appendix, and one that further extends
the first extension in section 2 above, is presented now. Here, instead of a structure of binary relations on the
arbitrary sets $X$ and $Y$ whose tensor product we define, we shall consider a more general structure given by {\it
generators}, as follows. \\

{\bf Definition 3.1.} \\

Given any set $X$, a mapping $\psi : {\cal P} ( X ) \longrightarrow {\cal P} ( X )$ will be
called a {\it generator}, if and only if \\

(3.1)~~~ $ \forall~~~ A \subseteq X ~:~ A \subseteq \psi ( A ) $ \\

and \\

(3.2)~~~ $ \forall~~~ A \subseteq A\,' \subseteq X ~:~ \psi ( A ) \subseteq \psi ( A\,' ) $ \\

{\bf Examples 3.1.} \\

1) A trivial example of generator is given by $\psi = id_{{\cal P} ( X )}$, that is, when $\psi
( A ) = A$, for $A \subseteq X$. \\

2) Let us now show that the concept of generator on a set $X$ is {\it more general} than that of a binary operation.
Indeed, given any binary operation $\alpha : X \times X \longrightarrow X$, we call a subset $A \subseteq X$ to be
$\alpha$-{\it stable}, if and only if \\

(3.3)~~~ $ x, y \in A ~~\Longrightarrow~~ \alpha ( x, y ) \in A $ \\

Obviously, $X$ is $\alpha$-stable, and the intersection of any family of $\alpha$-stable subsets is $\alpha$-stable.
Consequently, for every subset $A \subseteq X$, we can define the smallest $\alpha$-stable subset which contains it,
namely \\

(3.4)~~~ $ [ A ]_\alpha = \bigcap_{A \subseteq B,~ B ~\alpha-stable}~ B $ \\

Therefore, we can associate with $\alpha$ the mapping  $\psi_\alpha : {\cal P} ( X ) \longrightarrow {\cal P} ( X )$
defined by \\

(3.5)~~~ $ \psi_\alpha ( A ) = [ A ]_\alpha,~~~ A \subseteq X $ \\

which is obviously a generator. Furthermore, we have in view of (3.4) \\

(3.6)~~~ $ \forall~~~ A \subseteq X ~:~
       \psi_\alpha ( \psi_\alpha ( A ) ) = \psi_\alpha ( A ) $ \\

since as mentioned, $[ A ]_\alpha$ is $\alpha$-stable, and obviously $[ A ]_\alpha \subseteq [ A ]_\alpha$. \\

We note that, in general, the relation $\psi ( \psi ( A ) ) = \psi ( A )$, with $A \subseteq X$, need not hold for an
arbitrary generator $\psi$. \\

3) A particular case of 2) above is the following. Let $( S,\ast )$ be a semigroup with the
neutral element $e$. Then $[ \{ e \} ]_\ast = \{ e \}$, while for $a \in S,~ a \neq e$, we
have $[ \{ a \} ]_\ast = \{ a, a \ast a, a \ast a \ast a, \dots \}$. \\

For instance, if $( S, \ast ) = ( \mathbb{N}, + )$, then $[ \{ 0 \} ]_+ = \{ 0 \}$, while $[ \{ 1 \} ]_+ = \mathbb{N}
\setminus \{ 0 \} = \mathbb{N}_1$. \\

{\bf Definition 3.2.} \\

Given a generator $\psi$ on a set $X$. A binary operation $\alpha$ on $X$ is called {\it compatible} with $\psi$, if and
only if, see (3.5) \\

(3.7)~~~ $ \psi_\alpha ( A ) \subseteq \psi ( A ),~~~ A \subseteq X $

\hfill $\Box$ \\

Obviously, $\alpha$ is compatible with $\psi_\alpha$, for every binary operation $\alpha$ on $X$. \\

Let now be given two structures $( X, \psi )$ and $( Y, \varphi )$, where $\psi : {\cal P} ( X ) \longrightarrow
{\cal P} ( X ),~ \varphi : {\cal P} ( Y ) \longrightarrow {\cal P} ( Y )$ are arbitrary generators on $X$ and $Y$,
respectively. Let us again denote by $Z$ the set of all finite sequences of pairs \\

(3.8)~~~ $ ( x_1, y_1 ), \dots , ( x_n, y_n ) $ \\

where $n \geq 1$, while $x_i \in X,~ y_i \in Y$, with $1 \leq i \leq n$. Once more, we define on $Z$ the binary operation
$\gamma$  simply by the concatenation of the sequences (3.8). It follows that $\gamma$ is associative, therefore, each
sequence (3.8) can be written as \\

(3.9)~~~ $ ( x_1, y_1 ), \dots , ( x_n, y_n ) =
                 ( x_1, y_1 ) ~\gamma~ ( x_2, y_2 ) ~\gamma~ \ldots ~\gamma~ ( x_n, y_n ) $ \\

where for $n = 1$, the right hand term is understood to be simply $( x_1, y_1 )$. Obviously, if $X$ or $Y$ have at least
two elements, then $\gamma$ is not commutative. \\

Thus we have \\

(3.10)~~~ $ Z = \left \{ ( x_1, y_1 ) ~\gamma~ ( x_2, y_2 ) ~\gamma~ \ldots ~\gamma~ ( x_n, y_n ) ~~
                             \begin{array}{|l}
                              ~ n \geq 1 \\ \\
                              ~ x_i \in X,~ y_i \in Y,~ 1 \leq i \leq n
                             \end{array} \right \} $ \\ \\

which obviously gives \\

(3.11)~~~ $ X \times Y \subseteq Z $ \\

Now we define on $Z$ an equivalence relation $\approx_{\psi, \varphi}$ as follows. Two sequences in (3.8) are equivalent,
if and only if they are identical, or each can be obtained from the other by a finite number of applications of the
following operations \\

(3.12)~~~ permute pairs $( x_i, y_i )$ within the sequence \\ \\

(3.13)~~~ replace $( x_1, y_1 ) ~\gamma~ ( x\,'_1, y_1 ) ~\gamma~ ( x_2, y_2 ) ~\gamma~ \ldots ~\gamma~ ( x_n, y_n )$ \\
          \hspace*{1.5cm} with $( \alpha ( x_1, x\,'_1), y_1 ) ~\gamma~ ( x_2, y_2 ) ~\gamma~ \ldots , ( x_n, y_n )$,
          or vice-versa, \\
          \hspace*{1.5cm} where $\alpha$ is a binary operation on $X$ which is compatible
          \hspace*{1.5cm} with $\psi$ \\ \\

(3.14)~~~ replace $( x_1, y_1 ) ~\gamma~ ( x_1, y\,'_1 ) ~\gamma~ ( x_2, y_2 ) ~\gamma~ \ldots ~\gamma~ ( x_n, y_n )$ \\
          \hspace*{1.5cm} with $( x_1, \beta ( y_1, y\,'_1 ) ) ~\gamma~ ( x_2, y_2 ) ~\gamma~ \ldots ~\gamma~
          ( x_n, y_n )$, or vice-versa, \\
          \hspace*{1.5cm} where $\beta$ is a binary operation on $Y$ which is compatible
          \hspace*{1.5cm} with $\varphi$ \\ \\

We note again that, in view of Lemma 2.1., the binary relation $\approx_{\psi, \varphi}$ defined above is indeed an
equivalence relation on $Z$. \\

Finally, the {\it tensor product} of $( X, \psi )$ and $( Y, \varphi )$ is defined to be the quotient space \\

(3.15)~~~ $ X \bigotimes_{\psi, \varphi} Y = Z / \approx_{\psi, \varphi} $ \\

with the canonical quotient embedding, see (3.11) \\

(3.16)~~~ $ X \times Y \ni ( x, y ) \longmapsto x \bigotimes_{\psi, \varphi} y \in X \bigotimes_{\psi, \varphi} Y $ \\

where as in the usual case of tensor products, we denote by $x \bigotimes_{\psi, \varphi} y$, or simply  $x \bigotimes y$,
the equivalence class of $( x, y ) \in X \times Y \subseteq Z$. Again, similar with (2.9), it is easy to see that the mapping (3.16) is indeed
injective. \\

Obviously, the binary operation $\gamma$ on $Z$ will canonically lead by this quotient operation to a {\it commutative}
and {\it associative} binary operation on $X \bigotimes_{\psi, \varphi} Y$, which for convenience is denoted by the same
$\gamma$, although in view of (3.15), this time it depends on $\psi$ and $\varphi$, hence for the sake of rigor, it can be
written as $\gamma_{\psi, \varphi}$. \\

In this way, the elements of $X \bigotimes_{\alpha, \beta} Y$ are all the expressions \\

(3.17)~~~ $ x_1 \bigotimes_{\psi, \varphi} y_1 ~\gamma~ x_2 \bigotimes_{\psi, \varphi} y_2 ~\gamma~
                          \ldots ~\gamma~ x_n \bigotimes_{\psi, \varphi} y_n $ \\

with $n \geq 1$ and $x_i \in X,~ y_i \in Y$, for $1 \leq i \leq n$. \\

We conclude by \\

{\bf Theorem 3.1.} \\

The tensor products constructed in section 2 above are particular cases of those in this section. \\

{\bf Proof.} \\

Let be given two structures $( X, \alpha )$ and $( Y, \beta )$, where $\alpha : X \times X \longrightarrow X,~ \beta :
Y \times Y \longrightarrow Y$ are arbitrary binary operations on $X$ and $Y$, respectively. Then as in (3.5), we
associate with them the generators $\psi_\alpha$ and $\psi_\beta$ on $X$ and $Y$, respectively. \\

We show now that, for $z, z\,' \in Z$, we have the implication \\

(3.18)~~~ $ z \approx_{\alpha, \beta} z\,' ~~~\Longrightarrow~~~ z \approx_{\psi_\alpha, \psi_\beta} z\,' $ \\

Indeed, it is sufficient to prove that (2.6) implies (3.13), and (2.7) implies (3.14). And clearly, in both implications
we can  assume $n = 1$ without loss of generality. \\

Let us therefore be given $x, x\,' \in X,~ y \in Y$. If we assume (2.6), then we obtain \\

(3.19)~~~ $ ( \alpha ( x, x\,' ), y ) \approx_{\alpha, \beta} ( x, y ) ~\gamma~ ( x\,' y ) $ \\

But as we noted following Definition 3.2., $\alpha$ is compatible with $\phi_\alpha$, thus (3.19) gives \\

(3.20)~~~ $ ( \alpha ( x, x\,' ), y ) \approx_{\psi_\alpha, \psi_\beta} ( x, y ) ~\gamma~ ( x\,' y ) $ \\

and the implication (2.6) $\Longrightarrow$ (3.13) is proved. The proof of the implication (2.7) $\Longrightarrow$ (3.14)
is similar. \\ \\

{\bf 4. Universality Property of the Two Generalized Tensor \\
        \hspace*{0.45cm} Products} \\

A fundamental, and in fact, characterizing property of usual tensor products is their {\it universality} given in
(A1.4.3). Here we show that the two generalized tensor products introduced in sections 2 and 3 have similar universality
properties. \\

For that purpose, we have to introduce a few notions. Given three arbitrary sets $X, Y$ and $U$, as well as binary
operations $\alpha$ on $X$, $\beta$ on $Y$ and $\delta$ on $U$. A mapping $f : X \longrightarrow U$ is called an
$\alpha,\delta$-homomorphism, if and only if \\

(4.1)~~~ $ f ( \alpha ( x, x\,' ) ) = \delta ( f ( x ), f ( x\,' ) ),~~~ x, x\,' \in X $ \\

Further, a mapping $g : X \times Y \longrightarrow U$ is called an $\alpha,\beta,\delta$-homomorphism, if and only if \\

(4.2)~~~ $ g (  \alpha ( x, x\,' ), y ) = \delta ( g ( x, y ), g ( x\,', y ) ),~~~ x, x\,' \in X,~ y \in Y $ \\

and \\

(4.3)~~~ $ g ( x, \beta ( y, y\,' ) ) = \delta ( g ( x, y ), g ( x, y\,' ) ),~~~ x \in X,~ y, y\,' \in Y $ \\

Also, an $\alpha,\beta,\delta$-homomorphism $g : X \times Y \longrightarrow U$ is called an
$\alpha,\beta,\delta$-commuting homomorphism, if and only if \\

(4.4)~~~ $\delta ~\mbox{is commutative on}~ g ( X \times Y ) \subseteq U $ \\

{\bf Lemma 4.1.} \\

The canonical quotient embedding (2.9), namely \\

(4.5)~~~ $ X \times Y \ni ( x, y ) \stackrel{t_{\alpha, \beta}}\longmapsto
                     x \bigotimes_{\alpha, \beta} y \in X \bigotimes_{\alpha, \beta} Y $ \\

is an $\alpha, \beta, \gamma$-homomorphism. \\

{\bf Proof.} \\

Follows immediately from the definition of $\approx_{\alpha, \beta}$ in (2.5) - (2.7). \\

{\bf Theorem 4.1.} \\

The tensor products $X \bigotimes_{\alpha, \beta} Y$ have the following universality property :

(4.6) $\begin{array}{l}
              ~~~ \forall~~~ \alpha,\beta,\delta-\mbox{commuting homomorphism}~ g : X \times Y \longrightarrow U ~: \\ \\
              ~~~ \exists~ !~~ \gamma,\delta-\mbox{homomorphism}~ h_{\alpha, \beta} : X \bigotimes_{\alpha, \beta} Y
                            \longrightarrow U~ : \\ \\
              ~~~~~~~~~~ h_{\alpha, \beta} \circ t_{\alpha, \beta} ~=~ g
  \end{array} $ \\

where $t_{\alpha, \beta}$ is the canonical quotient embedding (2.9). In other words, the diagram commutes :

\begin{math}
\setlength{\unitlength}{1cm}
\thicklines
\begin{picture}(13,6)

\put(0,2.5){$(4.7)$}
\put(0.9,5){$X \times Y$}
\put(2.5,5.1){\vector(1,0){6.2}}
\put(9.2,5){$X \bigotimes_{\alpha, \beta} Y$}
\put(5,5.4){$t_{\alpha, \beta}$}
\put(1.7,4.5){\vector(1,-1){3.5}}
\put(9.5,4.5){\vector(-1,-1){3.5}}
\put(5.5,0.3){$U$}
\put(3,2.5){$g$}
\put(8.1,2.5){$\exists~!~~ h_{\alpha, \beta}$}

\end{picture}
\end{math} \\

{\bf Proof.} \\

We note that, for $x \in X,~ y \in Y$, the commutativity in (4.7) should give \\

(4.8)~~~ $ h_{\alpha, \beta} ( x \bigotimes_{\alpha, \beta} y ) = g ( x, y ) $ \\

which can therefore be taken as the definition of $h_{\alpha, \beta}$ on the range of the mapping $t_{\alpha, \beta}$.
However, in view of (2.10), this is sufficient in order to define $h_{\alpha, \beta}$ on the whole of $X
\bigotimes_{\alpha, \beta} Y$, simply by \\

(4.9)~~~ $ \begin{array}{l}
                h_{\alpha, \beta} ( x_1 \bigotimes y_1 ~\gamma~ x_2 \bigotimes y_2 ~\gamma~ \ldots ~\gamma~ x_n \bigotimes y_n ) = \\ \\
                ~~~~~~~~~~ = g ( x_1, y_1 ) ~\delta~ g ( x_2, y_2 ) ~\delta~ \ldots ~\delta~ g ( x_n, y_n )
           \end{array} $ \\

which will obviously make $h_{\alpha, \beta}$ into a $\gamma, \delta$-homomorphism. Let us note here that (4.9) and (2.5)
imply \\

(4.10)~~~ $ g ( x_1, y_1 ) ~\delta~ g ( x_2, y_2 ) =
                   g ( x_2, y_2 ) ~\delta~ g ( x_1, y_1 ),~~~ x_1, x_2 \in X,~ y_1, y_2 \in Y $ \\

which means that \\

(4.11)~~~ $ \delta ~\mbox{is commutative on}~ g ( X \times Y ) $ \\

this being the reason one has to ask that condition in (4.6), (4.7).

\hfill $\Box$ \\

Let us now turn to the universality property of the tensor products $X \bigotimes_{\psi, \varphi} Y$ defined in section 3. First, we note several useful
facts. Given any binary operation $\alpha$ on $X$ which is compatible with $\psi$, as well as any binary operation $\beta$ on $Y$ which is compatible
with $\varphi$, in view of (2.5) - (2.7) and (3.12) - (3.14), we obviously have for $z, z\,' \in Z$ \\

(4.12)~~~ $ z \approx_{\alpha, \beta} z\,' ~~~\Longrightarrow~~~  z \approx_{\psi, \varphi} z\,' $ \\

This clearly means that $X \bigotimes_{\alpha, \beta} Y = Z / \approx_{\alpha, \beta}$ is "larger" than $X
\bigotimes_{\psi, \varphi} Y = Z / \approx_{\alpha, \beta}$ in the sense that we have the surjective mapping \\

(4.13)~~~ $ X \bigotimes_{\alpha, \beta} Y \ni ( z )_{\approx_{\alpha, \beta}} \stackrel{i_{\alpha, \beta}}\longmapsto
                       ( z )_{\approx_{\alpha, \beta}} \in X \bigotimes_{\psi, \varphi} Y $ \\

where $( z )_{\approx_{\alpha, \beta}}$ denotes the $\approx_{\alpha, \beta}$ equivalence class of $z \in Z$, and
similarly with $( z )_{\approx_{\alpha, \beta}}$. Furthermore, $i_{\alpha, \beta}$ is a $\gamma, \gamma$-homomorphism,
more precisely, a $\gamma_{\alpha, \beta}, \gamma_{\psi, \varphi}$-homomorphism, see (2.10), (3.17). \\

Thus in view of (4.5), (3.16) and Theorem 4.1., we obtain the commutative diagram \\

\begin{math}
\setlength{\unitlength}{1cm}
\thicklines
\begin{picture}(15,6)

\put(0,2.5){$(4.14)$}
\put(0.9,5){$X \times Y$}
\put(2.5,5.1){\vector(1,0){8.1}}
\put(11,5){$X \bigotimes_{\psi, \varphi} Y$}
\put(6.1,5.4){$t_{\psi, \varphi}$}
\put(1.5,4.7){\vector(1,-1){4.7}}
%\put(11.7,4.7){\vector(-1,-1){4.7}}
\put(6.4,-0.45){$U$}
\put(3,2.3){$g$}
%\put(9.8,2.3){$???~ \exists~!~~ h_{\alpha, \beta}$}
\put(5.7,2.8){$X \bigotimes_{\alpha, \beta} Y$}
\put(2.2,4.7){\vector(2,-1){3.3}}
\put(3.8,4.1){$t_{\alpha, \beta}$}
\put(7.7,3.05){\vector(2,1){3.3}}
\put(8.8,4.1){$i_{\alpha, \beta}$}
\put(6.5,2.3){\vector(0,-1){2.2}}
\put(6.7,1.3){$h_{\alpha, \beta}$}

\end{picture}
\end{math} \\ \\

Let us note that in (4.14) the mappings $t_{\alpha, \beta}$ are $\alpha, \beta, \gamma$-homomorphisms, while $h_{\alpha,
\beta}$ are $\gamma, \delta$-homomorphisms. \\

Here we also note that the commutativity of the upper triangle in (4.14), namely

\begin{math}
\setlength{\unitlength}{1cm}
\thicklines
\begin{picture}(15,3)

\put(0,0.5){$(4.15)$}
\put(0.9,2){$X \times Y$}
\put(2.5,2.1){\vector(1,0){8.1}}
\put(11,2){$X \bigotimes_{\psi, \varphi} Y$}
\put(6.1,2.4){$t_{\psi, \varphi}$}
\put(5.7,-0.2){$X \bigotimes_{\alpha, \beta} Y$}
\put(2.2,1.7){\vector(2,-1){3.3}}
\put(3.8,1.1){$t_{\alpha, \beta}$}
\put(7.7,0.05){\vector(2,1){3.3}}
\put(8.8,1.1){$i_{\alpha, \beta}$}

\end{picture}
\end{math} \\ \\

can be construed as suggesting a {\it definition} of a $\psi, \varphi$-homomorphism, and in this case, specifically of the
mapping $t_{\psi, \varphi}$. More precisely, an arbitrary mapping \\

(4.16)~~~ $ s : X \times Y \longrightarrow X \bigotimes_{\psi, \varphi} Y $ \\

is called a $\psi, \varphi$-homomorphism, if and only if for ever binary operation $\alpha$ on $X$ compatible with $\psi$,
and every binary operation $\beta$ on $Y$ compatible with $\varphi$, we have a commutative diagram \\

\begin{math}
\setlength{\unitlength}{1cm}
\thicklines
\begin{picture}(15,3)

\put(0,0.5){$(4.17)$}
\put(0.9,2){$X \times Y$}
\put(2.5,2.1){\vector(1,0){8.1}}
\put(11,2){$X \bigotimes_{\psi, \varphi} Y$}
\put(6.1,2.4){$s$}
\put(5.7,-0.2){$X \bigotimes_{\alpha, \beta} Y$}
\put(2.2,1.7){\vector(2,-1){3.3}}
\put(3.8,1.1){$s_{\alpha, \beta}$}
\put(7.7,0.05){\vector(2,1){3.3}}
\put(8.8,1.1){$i_{\alpha, \beta}$}

\end{picture}
\end{math} \\ \\

where $s_{\alpha, \beta}$ is an $\alpha, \beta, \gamma_{\alpha, \beta}$-homomorphism. \\

We shall also need the following definition. An arbitrary mapping \\

(4.18)~~~ $ g : X \times Y \longrightarrow U $ \\

is called a $\psi, \varphi, \delta$-homomorphism, if and only if for ever binary operation $\alpha$ on $X$ compatible with
$\psi$, and every binary operation $\beta$ on $Y$ compatible with $\varphi$, we have \\

(4.17)~~~ $ h_{\alpha, \beta} \circ t_{\alpha, \beta} = g $ \\

where $t_{\alpha, \beta}$ is the canonical quotient embedding (2.9), while $h_{\alpha, \beta}$ corresponds to $g$
according to (4.7). \\

It follows that, for the universality property of the tensor products $X \bigotimes_{\psi, \varphi} Y$, all we need is to
complete the commutative diagram (4.14) with a mapping \\

(4.18)~~~ $ h_{\psi, \varphi} : X \bigotimes_{\psi, \varphi} Y \longrightarrow U $ \\

However, in view of the surjectivity of the mapping (4.13), such a completion is immediate. Indeed, we simply define \\

(4.19)~~~ $ h_{\psi, \varphi} ( i_{\alpha, \beta} ( x \bigotimes_{\alpha, \beta} y ) ) =
                 h_{\alpha, \beta} ( x \bigotimes_{\alpha, \beta} y ),~~~
                         x \bigotimes_{\alpha, \beta} y \in X \bigotimes_{\alpha, \beta} Y $ \\

and note that it is a correct definition. Furthermore, $h_{\psi, \varphi}$ turns out to be a $\gamma_{\psi, \varphi},
\delta$-homomorphism. \\

Thus we obtain \\

{\bf Theorem 4.2.} \\

The tensor products $X \bigotimes_{\psi, \varphi} Y$ have the following universality property :

(4.20) $ \begin{array}{l}
              ~~~ \forall~~~ \psi,\varphi,\delta-\mbox{commuting homomorphism}~ g : X \times Y \longrightarrow U~ : \\ \\
              ~~~ \exists~ !~~ \gamma_{\psi, \varphi},\delta-\mbox{homomorphism}~
                                   h_{\psi, \varphi} : X \bigotimes_{\psi, \varphi} Y \longrightarrow U~ : \\ \\
              ~~~~~~~~~~ h_{\psi, \varphi} \circ t_{\psi, \varphi} ~=~ g
          \end{array} $

\hfill $\Box$ \\

Here, similar with (4.4), a $\psi,\varphi,\delta$-homomorphism $g : X \times Y \longrightarrow U$ is called a
$\psi,\varphi,\delta$-commuting homomorphism, if and only if \\

(4.21)~~~ $ \delta ~\mbox{is commutative on}~ g ( X \times Y ) $ \\ \\

{\bf 5. Further Generalization} \\

The generalization of tensor products to the case of structures given by generators presented in section 3 above can easily be further extended. Indeed,
let $X,~ Y$ be arbitrary sets and let ${\cal A}$ be any set of binary operations on $X$, while correspondingly, ${\cal B}$ is any set of binary
operations on $Y$. \\

The constructions in (3.8) - (3.10) can again be implemented, since they only depend on the sets $X,~ Y$. \\

Now, we can define on $Z$ the equivalence relation $\approx_{{\cal A}, {\cal B}}$ as follows. Two sequences in (3.8) are equivalent, if and only if they
are identical, or each can be obtained from the other by a finite number of applications of the following operations \\

(5.1)~~~ permute pairs $( x_i, y_i )$ within the sequence \\ \\

(5.2)~~~ replace $( x_1, y_1 ) ~\gamma~ ( x\,'_1, y_1 ) ~\gamma~ ( x_2, y_2 ) ~\gamma~ \ldots ~\gamma~ ( x_n, y_n )$ \\
          \hspace*{1.5cm} with $( \alpha ( x_1, x\,'_1), y_1 ) ~\gamma~ ( x_2, y_2 ) ~\gamma~ \ldots , ( x_n, y_n )$,
          or vice-versa, \\
          \hspace*{1.5cm} where $\alpha \in {\cal A}$ \\ \\

(5.3)~~~ replace $( x_1, y_1 ) ~\gamma~ ( x_1, y\,'_1 ) ~\gamma~ ( x_2, y_2 ) ~\gamma~ \ldots ~\gamma~ ( x_n, y_n )$ \\
          \hspace*{1.5cm} with $( x_1, \beta ( y_1, y\,'_1 ) ) ~\gamma~ ( x_2, y_2 ) ~\gamma~ \ldots ~\gamma~
          ( x_n, y_n )$, or vice-versa, \\
          \hspace*{1.5cm} where $\beta \in {\cal B}$ \\ \\

Then once again, in view of Lemma 2.1., the binary relation $\approx_{{\cal A}, {\cal B}}$ defined above is indeed an equivalence relation on $Z$. \\

In rest, one can proceed with a construction similar with (3.15) - (3.17), and obtain the tensor product \\

(5.4)~~~ $ X \bigotimes_{{\cal A}, {\cal B}} Y = Z / \approx_{{\cal A}, {\cal B}} $ \\

with the canonical quotient embedding \\

(5.5)~~~ $ X \times Y \ni ( x, y ) \longmapsto x \bigotimes_{{\cal A}, {\cal B}} y \in X \bigotimes_{{\cal A}, {\cal B}} Y $ \\

and  $X \bigotimes_{{\cal A}, {\cal B}} Y$ being the set of all elements \\

(5.6)~~~ $ x_1 \bigotimes_{{\cal A}, {\cal B}} y_1 ~\gamma~ x_2 \bigotimes_{{\cal A}, {\cal B}} y_2 ~\gamma~
                                  \ldots ~\gamma~ x_n \bigotimes_{{\cal A}, {\cal B}} y_n $ \\

with $n \geq 1$ and $x_i \in X,~ y_i \in Y$, for $1 \leq i \leq n$. \\

Let us further note that, given sets ${\cal A} \subseteq {\cal A}\,'$ of binary operations on $X$, and correspondingly, sets ${\cal B} \subseteq
{\cal B}\,'$ of binary operations on $Y$, we have the surjective mapping \\

(5.7)~~~ $ X \bigotimes_{{\cal A}, {\cal B}} Y \ni ( z )_{{\cal A}, {\cal B}} \longmapsto
                         ( z )_{{\cal A}\,', {\cal B}\,'} \in  X \bigotimes_{{\cal A}\,', {\cal B}\,'} Y $ \\

where $( z )_{{\cal A}, {\cal B}}$ denotes the $\approx_{{\cal A}, {\cal B}}$ equivalence class of $z \in Z$, and similarly
with
$( z )_{{\cal A}\,', {\cal B}\,'}$. \\

Finally, similar with Theorem 4.2., and with the obvious definitions of the respective concepts of homomorphism, we have a
universality property of the tensor products (5.4), given in \\

{\bf Theorem 5.1.} \\

The tensor products $X \bigotimes_{{\cal A}, {\cal B}} Y$ have the following universality property :

(5.8) $ \begin{array}{l}
              ~~~ \forall~~~ \psi,\varphi,\delta-\mbox{commuting homomorphism}~ g : X \times Y \longrightarrow U~ : \\ \\
              ~~~ \exists~ !~~ \gamma_{{\cal A}, {\cal B}},\delta-\mbox{homomorphism}~
                                   h_{{\cal A}, {\cal B}} : X \bigotimes_{{\cal A}, {\cal B}} Y \longrightarrow U~ : \\ \\
              ~~~~~~~~~~ h_{{\cal A}, {\cal B}} \circ t_{{\cal A}, {\cal B}} ~=~ g
          \end{array} $ \\ \\

{\bf Appendix. Definition of Usual Tensor Products \\
               \hspace*{2.25cm} of Vector Spaces} \\

For convenience, we recall here certain main features of the usual tensor product of vector
spaces, and relate them to certain properties of Cartesian products. \\

Let $\mathbb{K}$ be a field and $E,F, G$ vector spaces over $\mathbb{K}$. \\

{\bf A1.1. Cartesian Product of Vector Spaces} \\

Then $E \times F$ is the vector space over $\mathbb{K}$ where the operations are given by \\

$~~~~~~ \lambda ( x, y ) + \mu ( u, v ) ~=~ ( \lambda x + \mu u, \lambda y + \mu v ) $ \\

for any $x, y \in E,~ u, v \in F,~ \lambda, \mu \in \mathbb{K}$. \\ \\

{\bf A1.2. Linear Mappings} \\

Let ${\cal L} ( E, F )$ be the set of all mappings \\

$~~~~~~ f : E ~\longrightarrow~ F $ \\

such that \\

$~~~~~~ f ( \lambda x + \mu u ) ~=~ \lambda f ( x ) + \mu f ( u ) $ \\

for $u, v \in E,~ \lambda, \mu \in \mathbb{K}$. \\ \\

{\bf A1.3. Bilinear Mappings} \\

Let ${\cal L} ( E, F; G )$ be the set of all mappings \\

$~~~~~~ g : E \times F ~\longrightarrow~ G $ \\

such that for $x \in E$ fixed, the mapping $F \ni y \longmapsto g ( x, y ) \in G$ is linear in
$y$, and similarly, for $y \in F$ fixed, the mapping $E \ni x \longmapsto g ( x, y ) \in G$ is
linear in $x \in E$. \\

It is easy to see that \\

$~~~~~~ {\cal L} ( E, F; G ) ~=~ {\cal L} ( E, {\cal L} ( F, G ) ) $ \\ \\

{\bf A1.4. Tensor Products} \\

The aim of the tensor product $E \bigotimes F$  is to establish a close connection between the
{\it bilinear} mappings in ${\cal L} ( E, F; G )$ and the {\it linear} mappings in
${\cal L} ( E \bigotimes F , G )$. \\

Namely, the {\it tensor product} $E \bigotimes F$ is : \\

(A1.4.1)~~~ a vector space over $\mathbb{K}$, together with \\

(A1.4.2)~~~ a bilinear mapping $t : E \times F ~\longrightarrow~ E \bigotimes F$, such that we \\
          \hspace*{1.6cm} have the following : \\ \\

{\bf UNIVERSALITY PROPERTY} \\

$\begin{array}{l}
    ~~~~~~~~~~ \forall~~~ V ~\mbox{vector space over}~ \mathbb{K},~~
                   g \in {\cal L} ( E, F; V ) ~\mbox{bilinear mapping}~ : \\ \\
    ~~~~~~~~~~ \exists~ !~~ h \in {\cal L} ( E \bigotimes F, V )
                                         ~\mbox{linear mapping}~ : \\ \\
    ~~~~~~~~~~~~~~~~ h \circ t ~=~ g
  \end{array} $ \\ \\

or in other words : \\

(A1.4.3)~~~ the diagram commutes

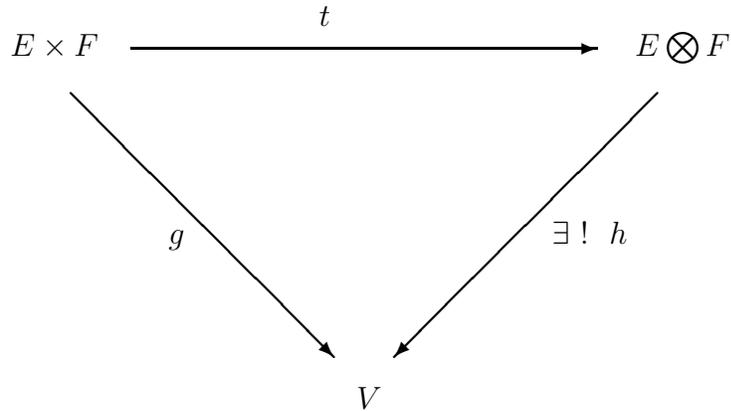
\begin{math}
\setlength{\unitlength}{1cm}
\thicklines
\begin{picture}(13,7)

\put(0.9,5){$E \times F$}
\put(2.5,5.1){\vector(1,0){6.2}}
\put(9.2,5){$E \bigotimes F$}
\put(5,5.4){$t$}
\put(1.7,4.5){\vector(1,-1){3.5}}
%\psline[linestyle=dashed]{->}(9.5,4.5)(6,1)
\put(9.5,4.5){\vector(-1,-1){3.5}}
\put(5.5,0.3){$V$}
\put(3,2.5){$g$}
\put(8.1,2.5){$\exists~!~~ h$}

\end{picture}
\end{math}

and \\

(A1.4.4)~~~ the tensor product $E \bigotimes F$ is {\it unique} up to vector \\
          \hspace*{1.6cm} space isomorphism. \\

Therefore we have the {\it injective} mapping \\

$~~~~~~ {\cal L} ( E, F; V ) \ni g ~\longmapsto~ h \in {\cal L} ( E \bigotimes F, V )
                                                    ~~~~\mbox{with}~~~~ h \circ t ~=~ g $ \\

The converse mapping \\

$~~~~~~ {\cal L} ( E \bigotimes F, V ) \ni h ~\longmapsto~
                         g ~=~ h \circ t \in {\cal L} ( E, F; V ) $ \\

obviously exists. Thus we have the {\it bijective} mapping \\

$~~~~~~ {\cal L} ( E \bigotimes F, V ) \ni h ~\longmapsto~
                         g ~=~ h \circ t \in {\cal L} ( E, F; V ) $ \\ \\

{\bf A1.5. Lack of Interest in ${\cal L} ( E \times F, G )$} \\

Let $f \in {\cal L} ( E \times F, G )$ and $( x, y ) \in E \times F$, then $( x, y ) = ( x, 0 )
+ ( 0, y )$, hence \\

$~~~~~~ f ( x, y ) ~=~ f ( ( x, 0 ) + ( 0, y ) ) ~=~ f ( x, 0 ) ~+~ f ( 0, y ) $ \\

thus $f ( x, y )$ depends on $x$ and $y$ in a {\it particular} manner, that is, separately on
$x$, and separately on $y$. \\ \\

{\bf A1.6. Universality Property of Cartesian Products} \\

Let $X, Y$ be two nonvoid sets. Their cartesian product is : \\

(A1.6.1)~~~ a set $X \times Y$, together with \\

(A1.6.2)~~~ two projection mappings~ $p_X : X \times X ~\longrightarrow~ X, \\
          \hspace*{1.6cm} p_Y : X \times Y ~\longrightarrow~ Y$, such that we have the
          following : \\ \\

{\bf UNIVERSALITY PROPERTY} \\

$\begin{array}{l}
    ~~~~~~~~~~ \forall~~~ Z ~\mbox{nonvoid set},~~ f : Z ~\longrightarrow~ X,~~
                          g : Z ~\longrightarrow~ Y~ : \\ \\
    ~~~~~~~~~~ \exists~ !~~ h  : Z ~\longrightarrow~ X \times Y~ : \\ \\
    ~~~~~~~~~~~~~~~~ f ~=~ p_X \circ h,~~~ g ~=~ p_Y \circ h
  \end{array} $ \\ \\

or in other words : \\

(A1.6.3)~~~ the diagram commutes

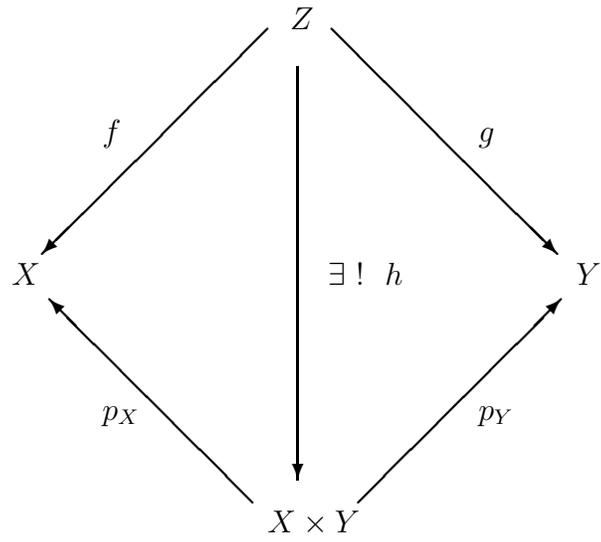
\begin{math}
\setlength{\unitlength}{1cm}
\thicklines
\begin{picture}(20,9)

\put(6.5,3.6){$\exists~ ! ~~ h$}
\put(6.1,6.5){\vector(0,-1){5.5}}
\put(6,7){$Z$}
\put(3.5,5.5){$f$}
\put(5.7,7){\vector(-1,-1){3}}
\put(8.5,5.5){$g$}
\put(6.55,7){\vector(1,-1){3}}
\put(2.3,3.6){$X$}
\put(9.8,3.6){$Y$}
\put(3.5,1.8){$p_X$}
\put(5.5,0.7){\vector(-1,1){2.7}}
\put(8.5,1.8){$p_Y$}
\put(6.9,0.7){\vector(1,1){2.7}}
\put(5.7,0.3){$X \times Y$}

\end{picture}
\end{math} \\ \\

{\bf A1.7. Cartesian and Tensor Products seen together} \\

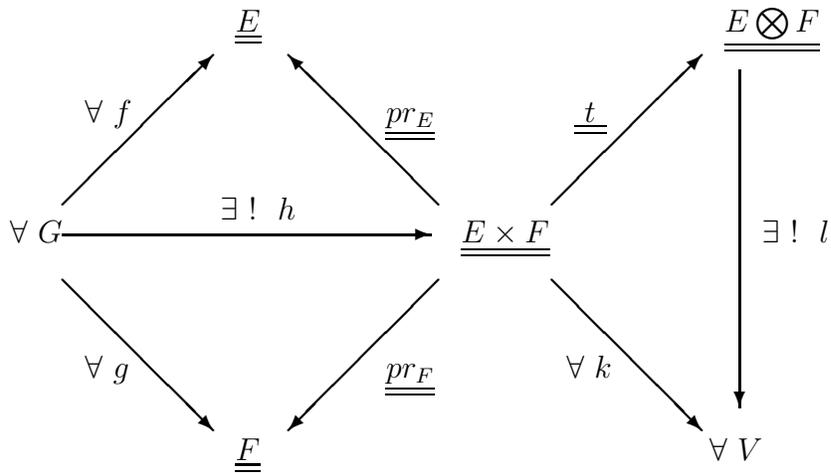
\begin{math}
\setlength{\unitlength}{1cm}
\thicklines
\begin{picture}(30,10)

\put(-0.2,6){$\forall~G$}
\put(0.5,6.5){\vector(1,1){2}}
\put(0.8,7.6){$\forall~f$}
\put(2.8,8.8){$\underline{\underline{E}}$}
\put(0.5,5.5){\vector(1,-1){2}}
\put(0.8,4.2){$\forall~g$}
\put(2.8,3.1){$\underline{\underline{F}}$}
\put(5.8,6){$\underline{\underline{E \times F}}$}
\put(5.5,6.5){\vector(-1,1){2}}
\put(4.8,7.6){$\underline{\underline{pr_E}}$}
\put(5.5,5.5){\vector(-1,-1){2}}
\put(4.8,4.2){$\underline{\underline{pr_F}}$}
\put(0.5,6.1){\vector(1,0){4.9}}
\put(2.6,6.3){$\exists~ !~~ h$}
\put(7,6.5){\vector(1,1){2}}
\put(7.3,7.6){$\underline{\underline{~t~}}$}
\put(7,5.5){\vector(1,-1){2}}
\put(7.2,4.2){$\forall~k$}
\put(9.3,8.8){$\underline{\underline{E \bigotimes F}}$}
\put(9.1,3.1){$\forall~V$}
\put(9.5,8.3){\vector(0,-1){4.5}}
\put(9.8,6){$\exists~ ! ~~ l$}

\end{picture}
\end{math} \\ \\

{\bf Acknowledgement} \\

Many thanks to my colleague Gusti van Zyl for pointing out the need for the condition of {\it commutative homeomorphism}
in the universality property of tensor products. \\ \\

{~} \\ \\ \\

{\large \bf PART 2 : Superpositions} \\ \\

{\bf 1. Superposition as a Particular case of Extended Tensor \\
        \hspace*{0.45cm} Products} \\

Recently in [10-13] the concepts of {\it tensor products} and {\it
entanglement} were extended far beyond quantum systems, being
defined for rather arbitrary systems which need not have any kind of
quantum nature. Here it is shown that {\it superposition} is in fact
a particular case of the mentioned
extended concept of tensor product, therefore, it can be defined far beyond linear systems for which it is usually considered. \\

{\bf Extended Tensor Products} \\

For convenience, let us recall the extended definition of
entanglement. Let be given an arbitrary family of sets $( X_i )_{i
\in I}$, that is,
both the index set $I$ and the sets $X_i$, with $i \in I$, are arbitrary. \\

Further, let ${\cal A}$ be a family of $(n+m)$-ary relations $A
\subseteq ( \prod_{i \in I} X_i )^n \times ( \prod_{i \in I} X_i
)^m$. Here $n, m \geq
0,~ n + m \geq 1$, may depend on $A \in {\cal A}$. \\

For clarity, we note that, given an $(n+m)$-ary relation $A \in {\cal A}$, its elements are sequences of the form \\

(1.1)~~~ $ ( x_{1,i} )_{i \in I}, \ldots , ( x_{n,i} )_{i \in I}, ( x\,'_{1,i} )_{i \in I}, \ldots , ( x\,'_{m,i} )_{i \in I} $ \\

where $x_{\nu,i}, x\,'_{\mu, i} \in X_i$, with $1 \leq \nu \leq n,~ 1 \leq \mu \leq m,~ i \in I$. \\

We start the construction of the extended tensor product, and at
that stage, we do so based alone on $\prod_{i \in I} X_i$. Namely,
we consider the
semigroup \\

(1.2)~~~ $ ( Z, \gamma ) $ \\

given by the {\it free semigroup}, see Appendix, generated by $\prod_{i \in I} X_i$. In other words, $Z$ has as elements all the finite sequences \\

(1.3)~~~ $ ( x_{1,i} )_{i \in I}, \ldots , ( x_{h,i} )_{i \in I} $ \\

where $h \geq 1$, while $x_{k, i} \in X_i$, with $1 \leq k \leq h,~
i \in I$. Further, the binary operation $\gamma$ on $Z$ is simply
the concatenation
of the sequences (1.3). It follows that $\gamma$ is {\it associative}, therefore, each  sequence (1.3) can be written as \\

(1.4)~~~ $ ( x_{1,i} )_{i \in I} ~\gamma~ \ldots ~\gamma~ ( x_{h,i} )_{i \in I} $ \\

where for $h = 1$, the expression is understood to be simply $(
x_{1,i} )_{i \in I}$. Obviously, if one of the $X_i$, with $i \in
I$, has at least two
elements, then $\gamma$ is {\it not} commutative. \\

Thus we have \\

(1.5)~~~ $ Z = \{~ ( x_{1,i} )_{i \in I} ~\gamma~ \ldots ~\gamma~ ( x_{h,i} )_{i \in I} ~~|~~ h \geq 1 ~\} $ \\

which clearly gives the {\it injective} mapping \\

(1.6)~~~ $ \prod_{i \in I} X_i \ni ( x_i )_{i \in I} \longmapsto ( x_i )_{i \in I} \in Z $ \\

Now, with the help of the family ${\cal A}$ of multi-arity
relations, we can define on $Z$ the {\it equivalence} relation
$\approx_{\cal A}$ as follows. Two sequences in (1.3) - (1.5) are
equivalent, if and only if they are identical, or each can be
obtained from the other by a finite number of
applications of the following operations : \\

(1.7)~~~ permute two $( x_{k,i} )_{i \in I},~ ( x_{l,i} )_{i \in I}$ within the sequence \\

(1.8)~~~ replace $( x_{1,i} )_{i \in I} ~\gamma~ \ldots ~\gamma~ (
x_{n,i} )_{i \in I} ~\gamma ~( y_i )_{i \in I}
           ~\gamma~ \ldots $

           ~~\hspace*{1cm} with $( x\,'_{1,i} )_{i \in I} ~\gamma~ \ldots ~\gamma~ ( x\,'_{m,i} )_{i \in I} ~\gamma~ ( y_i )_{i \in I} ~\gamma~ \ldots $

           ~~\hspace*{1cm} or vice-versa, where

           ~~\hspace*{1cm} $( x_{1,i} )_{i \in I}, \ldots , ( x_{n,i} )_{i \in I}, \ldots ,
           ( x\,'_{1,i} )_{i \in I}, \ldots , ( x\,'_{m,i} )_{i \in I} \in A $,

           ~~\hspace*{1cm} for some $A \in {\cal A}$ \\

As is well known, see Appendix, in this way one indeed obtains on
$Z$ an equivalence relation $\approx_{\cal A}$, by the general set
theoretic method
called {\it transitive closure}. \\

Thus one can define the extended {\it tensor product} by the usual set theoretic quotient construction, see Appendix, namely \\

{\bf Definition 1.1.} \\

Given an arbitrary family of sets $( X_i )_{i \in I}$, and a family
${\cal A}$ of $(n+m)$-ary relations $A \subseteq ( \prod_{i \in I}
X_i )^n \times (
\prod_{i \in I} X_i )^m$, where $n, m \geq 0,~ n + m \geq 1$, may depend on $A \in {\cal A}$. Then we define the extended {\it tensor product} \\

(1.9)~~~ $ \bigotimes_{{\cal A}, i \in I} X_i = Z / \approx_{\cal A} $ \\

with the resulting natural mapping \\

(1.10)~~~ $ \prod_{i \in I} X_i \ni ( x_i )_{i \in I}
\stackrel{\iota}\longmapsto
                      \bigotimes_{{\cal A}, i \in I} x_i \in \bigotimes_{{\cal A}, i \in I} X_i $ \\

where $\bigotimes_{{\cal A}, i \in I} x_i$ is the coset of all
elements in $Z$ which are equivalent with $( x_i )_{i \in I} \in Z$,
in the sense of the equivalence relation $\approx_{\cal A}$.

\hfill $\Box$ \\

Related to this natural mapping $\iota$ in (1.10) and  as an obvious consequence of (1.8), we have \\

{\bf Lemma 1.1.} \\

If $m, n \geq 2$ for every $(n+m)$-ary relation $A \subseteq (
\prod_{i \in I} X_i )^n \times ( \prod_{i \in I} X_i )^m$ in ${\cal
A}$, then the natural mapping $\iota$ in (1.10) is {\it injective},
thus it is an {\it embedding} of the cartesian product $\prod_{i \in
I} X_i$ into the extended tensor product $\bigotimes_{{\cal A}, i
\in I} X_i$.

\hfill $\Box$ \\

However, there are as well other cases when the natural mapping $\iota$ in (1.10) is injective, [10-13]. \\

An important feature of the equivalence relation $\approx_{\cal A}$
is that it is a {\it congruence} on the semigroup $( Z, \gamma )$,
see Appendix, a property which follows immediately from (1.7),
(1.8). Consequently, the above quotient construction in (1.9) leads
naturally to a {\it commutative
semigroup} structure on the extended tensor product $\bigotimes_{{\cal A}, i \in I} X_i$, namely, we are led to, see (1.2) \\

{\bf Definition 1.2.} \\

Given an arbitrary family of sets $( X_i )_{i \in I}$, and a family
${\cal A}$ of $(n+m)$-ary relations $A \subseteq ( \prod_{i \in I}
X_i )^n \times ( \prod_{i \in I} X_i )^m$, where $n, m \geq 0,~ n +
m \geq 1$, may depend on $A \in {\cal A}$. Then we define the {\it
algebraic} structure on the
extended {\it tensor product} by \\

(1.11)~~~ $ ( \bigotimes_{{\cal A}, i \in I} X_i, \gamma_{\cal A} ) = ( Z, \gamma ) / \approx_{\cal A} $ \\

with the {\it associative} and {\it commutative} binary operation $\gamma_{\cal A}$ on $\bigotimes_{{\cal A}, i \in I} X_i$ given simply by \\

(1.12)~~~ $ ( z )_{\cal A} ~\gamma_{\cal A}~ ( z\,' )_{\cal A} = ( z ~\gamma~ z\,')_{\cal A},~~~ z, z\,' \in Z $ \\

where $( z )_{\cal A}$ denotes the $\approx_{\cal A}$ equivalence
class of $z \in Z$, and similarly with $( z\,' )_{\cal A}$ and $( z
~\gamma~ z\,')_{\cal A}$.

\hfill $\Box$ \\

Consequently, and according to the general properties of quotient
constructions, see Appendix, we obtain the extended tensor product
$\bigotimes_{{\cal A}, i \in I} X_i$ as the  set of all elements \\

(1.13)~~~ $ ( \bigotimes_{{\cal A}, i \in I} x_{1,i} ) ~\gamma_{\cal
A}~ \ldots ~\gamma_{\cal A}~
                       ( \bigotimes_{{\cal A}, i \in I} x_{h,i} ) $ \\

with $h \geq 1$ and $x_{k, i} \in X_i$, for $1 \leq k \leq h,~ i \in I$. \\

Usually, we shall simply write $\gamma$ instead of $\gamma_{\cal A}$, and the notation in (1.13) will further be simplified to \\

(1.14)~~~ $ ( x_{1, 1} \bigotimes_{\cal A} \ldots \bigotimes_{\cal
A} x_{1, L} ) ~\gamma~ \ldots ~\gamma~
                ( x_{h, 1} \bigotimes_{\cal A} \ldots \bigotimes_{\cal A} x_{h, L} ) $ \\

Let us further note that, given sets ${\cal A} \subseteq {\cal A}\,'$ of multi-arity relations as above, we have the surjective mapping \\

(1.15)~~~ $ \bigotimes_{{\cal A}, i \in I} X_i \ni ( z )_{\cal A}
\longmapsto ( z )_{{\cal A}\,'} \in
              \bigotimes_{{\cal A}\,', i \in I} X_i $ \\

where $( z )_{\cal A}$ denotes the $\approx_{\cal A}$ equivalence class of $z \in Z$, and similarly with $( z )_{{\cal A}\,'}$. \\

{\bf Extended Entanglement} \\

The extended concept of {\it entanglement} in the context of the above extended concept of tensor product, is given by \\

{\bf Definition 1.3.} \\

The elements in \\

(1.16)~~~ $ \bigotimes_{{\cal A}, i \in I} X_i ~\setminus~ \iota ( \prod_{i \in I} X_i ) $ \\

are called the {\it entangled} elements of the tensor product
$\bigotimes_{{\cal A}, i \in I} X_i$.

\hfill $\Box$ \\

Clearly, such entangled elements exist, if and only if \\

(1.17)~~~ $ \bigotimes_{{\cal A}, i \in I} X_i ~\setminus~ \iota ( \prod_{i \in I} X_i ) \neq \phi $ \\

and in order that (1.17) hold, it is {\it not} necessary for the mapping $\iota$ in (1.10) to be injective. \\

{\bf Extended Superposition} \\

Now we can turn to the extended concept of {\it superposition}.
Obviously, superposition is supposed to take place within one single
system. Consequently, we consider the {\it particular} case of the
extended tensor product (1.9), (1.11), when the index set $I$ has
one single element, thus the
family of sets $( X_i )_{i \in I}$ reduces to one single arbitrary set $X$ and we then have \\

(1.18)~~~ $ \bigotimes_{\cal A} X $ \\

Here we note, however, that even in such a particular case of
extended tensor product, the family ${\cal A}$ of $(n+m)$-ary
relations $A \subseteq X^n
\times X^m$ can nevertheless be rather arbitrary. Therefore, there exists the possibility for a large variety of superpositions on a given set $X$. \\

Let us see now what becomes of (1.5) - (1.14) in such a case. It is easy to see that we shall have \\

(1.19)~~~ $ Z = \{~ x_1 ~\gamma~ \ldots ~\gamma~ x_h ~~|~~ h \geq 1,~~ x_1, \ldots , x_h \in X ~\} $ \\

which clearly gives the {\it injective} mapping \\

(1.20)~~~ $ X \ni x \longmapsto x \in Z $ \\

Further, {\it equivalence} relation $\approx_{\cal A}$ on $Z$ is
defined as follows. Two sequences in (1.19) are equivalent, if and
only if they are
identical, or each can be obtained from the other by a finite number of applications of the following operations : \\

(1.21)~~~ permute two $x_k,~ x_l$ within the sequence \\

(1.22)~~~ replace $x_1 ~\gamma~ \ldots ~\gamma~ x_n ~\gamma ~ y
~\gamma~ \ldots $

           ~~\hspace*{1.3cm} with $x\,'_1 ~\gamma~ \ldots ~\gamma~ x\,'_m ~\gamma~ y ~\gamma~ \ldots $

           ~~\hspace*{1.3cm} or vice-versa, where

           ~~\hspace*{1.3cm} $( x_1, \ldots , x_n, x\,'_1, \ldots , x\,'_m ) \in A $,

           ~~\hspace*{1.3cm} for some $A \in {\cal A}$ \\

Thus (1.18) follows by the usual set theoretic quotient construction, see Appendix \\

(1.23)~~~ $ \bigotimes_{\cal A} X = Z / \approx_{\cal A} $ \\

with the resulting natural mapping \\

(1.24)~~~ $ X \ni x \stackrel{\iota}\longmapsto \bigotimes_{\cal A} x \in \bigotimes_{\cal A} X $ \\

where $\bigotimes_{\cal A} x$ is the coset of all elements in $Z$
which are equivalent with $x \in Z$, in the sense of the equivalence
relation
$\approx_{\cal A}$. \\

Also, the equivalence relation $\approx_{\cal A}$ is a {\it
congruence} on the semigroup $( Z, \gamma )$, see Appendix, a
property which follows immediately from (1.21), (1.22).
Consequently, the above quotient construction in (1.23) leads
naturally to a {\it commutative semigroup} structure on
the tensor product $\bigotimes_{\cal A} X$, namely, we obtain  \\

(1.25)~~~ $ ( \bigotimes_{\cal A} X \gamma_{\cal A} ) = ( Z, \gamma ) / \approx_{\cal A} $ \\

with the {\it associative} and {\it commutative} binary operation $\gamma_{\cal A}$ on $\bigotimes_{\cal A} X$ given simply by \\

(1.26)~~~ $ ( z )_{\cal A} ~\gamma_{\cal A}~ ( z\,' )_{\cal A} = ( z ~\gamma~ z\,')_{\cal A},~~~ z, z\,' \in Z $ \\

where $( z )_{\cal A}$ denotes the $\approx_{\cal A}$ equivalence
class of $z \in Z$, and similarly with $( z\,' )_{\cal A}$ and $( z
~\gamma~
z\,')_{\cal A}$. \\

And thus according to the general properties of quotient
constructions, see Appendix, we obtain the tensor product
$\bigotimes_{\cal A} X$ as the set
of all elements \\

(1.27)~~~ $ ( \bigotimes_{\cal A} x_1 ) ~\gamma_{\cal A}~ \ldots
~\gamma_{\cal A}~
                       ( \bigotimes_{\cal A} x_h ) $ \\

with $h \geq 1$ and $x_k \in X$, for $1 \leq k \leq h$. Usually, we
shall simply write $\gamma$ instead of $\gamma_{\cal A}$, and the
notation in (1.27)
will further be simplified to \\

(1.28)~~~ $ ( \bigotimes_{\cal A} x_1 ) ~\gamma~ \ldots ~\gamma~ ( \bigotimes_{\cal A} x_h ) $ \\

with $h \geq 1$ and $x_k \in X$, for $1 \leq k \leq h$. \\ \\

{\bf 2. Examples of Classical Type Superpositions} \\

{\bf The Case of Arbitrary Binary Operations} \\

Let us consider the case when the family ${\cal A}$ of $(n+m)$-ary
relations $A \subseteq X^n \times X^m$ consists of one single such
$(2+1)$-ary
relation $A \subseteq X^2 \times X$, defined as follows. Let $\alpha : X \times X \longrightarrow X$ be any binary operation on $X$, and take \\

(2.1)~~~ $ A = \{~ ( x, y, z ) \in X^3 ~~|~~ \alpha ( x, y ) = z ~\} $ \\

Then (1.22) obviously becomes \\

(2.2)~~~ replace $x ~\gamma~ y ~\gamma ~ u ~\gamma~ \ldots $

           ~~\hspace*{1cm} with $\alpha ( x, y ) ~\gamma~ u ~\gamma~ \ldots $

           ~~\hspace*{1cm} or vice-versa \\

which means that, in terms of (1.28), we have in the tensor product $\bigotimes_{\cal A} X$ the relations \\

(2.3)~~~ $ ( \bigotimes_{\cal A} x ) ~\gamma~ ( \bigotimes_{\cal A}
y ) ~\gamma~ ( \bigotimes_{\cal A} u ) ~\gamma~ \ldots ~~=~~
              ( \bigotimes_{\cal A} \alpha ( x, y ) ) ~\gamma~ ( \bigotimes_{\cal A} u ) ~\gamma~ \ldots $ \\

for $x, y, u \in X$. In particular, we have in the tensor product $\bigotimes_{\cal A} X$ the relations \\

(2.4)~~~ $ ( \bigotimes_{\cal A} x ) ~\gamma~ ( \bigotimes_{\cal A} y ) ~~=~~ \bigotimes_{\cal A} \alpha ( x, y ),~~~ x, y \in X $ \\

In order to further clarify the effect of (2.2), we return to
(1.23), (1.24) and ask the following two questions

\begin{itemize}

\item when is the mapping $\iota$ in (1.24) {\it injective} ?

\item in such a case, how much $\bigotimes_{\cal A} X$ is larger than $\iota ( X )$ ?

\end{itemize}

The answer to these two questions is given by \\

{\bf Theorem 2.1.} \\

The mapping \\

(2.5)~~~ $ X \ni x \stackrel{\iota}\longmapsto \bigotimes_{\cal A} x \in \bigotimes_{\cal A} X $ \\

is always {\it surjective}, that is, we always have \\

(2.6)~~~ $ \iota \,( X ) = \bigotimes_{\cal A} X $ \\

Furthermore \\

(2.7)~~~ $ \iota $ is {\it injective} $~~~\Longleftrightarrow~~~ \iota $ is {\it bijective}  $~~~\Longleftrightarrow~~~ \alpha $ is associative \\

{\bf Proof} \\

The condition $\iota \,( X ) = \bigotimes_{\cal A} X$, is obviously equivalent with \\

$~~~~~~ \forall~ z \in Z ~:~ \exists~ x \in X ~:~ x ~\approx_{\cal A}~ z $ \\

thus further equivalent with \\

$~~~~~ \forall~ x_1, \ldots , x_h \in X ~:~ \exists~ x \in X ~:~ x ~\approx_{\cal A}~ x_1 ~\gamma~ \ldots ~\gamma~ x_h $ \\

which in view of (2.2) obviously holds whether or not $\alpha$ is associative, since for $h \geq 2$, we can take for instance \\

$~~~~~~ x = \alpha ( x_1, \alpha ( x_2, \ldots \alpha ( x_{h-1}, x_h ) \ldots ) $ \\

Regarding (2.7), it is easy to note that the associativity of
$\alpha$ is a necessary condition for the injectivity of $\iota$.
Indeed, let $a, b, c \in
X$, such that \\

$~~~~~~ u = \alpha ( a, \alpha ( b, c ) ) \neq \alpha ( \alpha ( a, b ), c ) = v $ \\

then in view of (2.2), we have \\

$~~~~~~ u ~\approx_{\cal A}~ a ~\gamma~ \alpha ( b, c )
~\approx_{\cal A}~ a ~\gamma~ b ~\gamma~ c ~\approx_{\cal A}~ \alpha
( a, b ) ~\gamma~ c
                                          ~\approx_{\cal A}~ v $ \\

Conversely, let us assume that $\alpha$ is associative. Then, due to the associativity of $\gamma$ itself, we obviously have the implication \\

$~~~~~~ \bigotimes_{\cal A} u ~=~ \bigotimes_{\cal A} v ~~~\Longrightarrow~~~ u ~=~ v $ \\

for $u, v \in X$. \\

{\bf The Particular Case of Usual Superpositions} \\

Let us illustrate the above by showing that the usual superposition
can indeed be obtained as a {\it particular} case. For that purpose,
let $( S, \ast )$ be a semigroup. In this case we can take $X = S$
and $\alpha = \ast : S \times S \longrightarrow S$ which is now
associative. Thus in view of
Theorem 2.1., it is easy to see that \\

(2.8)~~~ $ \begin{array}{l}
                 \mbox{the mapping}~~ \iota  ~~\mbox{in (2.5) is a {\it semigroup isomorphism}} \\
                 \mbox{between}~~ ( S, \ast ) ~~\mbox{and}~~ ( \bigotimes_{\cal A} X, \gamma )
           \end{array} $ \\

Indeed, let $x, y \in X$, then in view of (1.24), (2.4), we have
$\iota ( \alpha ( x, y ) ) = \bigotimes_{\cal A} \alpha ( x, y ) = (
\bigotimes_{\cal
A} x ) ~\gamma~ ( \bigotimes_{\cal A} y )$. \\

{\bf The Case of Non-Associative Binary Operations} \\

In view of (2.7), there is an interest in binary operations $\alpha$
which are {\it not} associative, thus for which $\bigotimes_{\cal A}
X$ may possibly have a {\it smaller} cardinal than $X$. Let us
therefore consider such a class of examples, in order to gain some
insight into the respective
situation. \\

We take $X = \mathbb{N}$ and $\alpha : X \times X \longrightarrow X$ defined by \\

(2.9)~~~ $ \alpha ( x, y ) = a x + b y,~~~ x, y \in \mathbb{N} $ \\

where $a, b \in \mathbb{N},~ a, b \geq 1,~ a + b \geq 3$ are given. Thus $\alpha$ is clearly non-associative. Then (2.2) becomes \\

(2.10)~~~ replace $x ~\gamma~ y ~\gamma ~ u ~\gamma~ \ldots $

           ~~~\hspace*{1cm} with $( a x + b y  ) ~\gamma~ u ~\gamma~ \ldots $

           ~~~\hspace*{1cm} or vice-versa \\

Then for $x_1, \ldots , x_h \in X$, with $h \geq 2$, we have by successively grouping the terms in $x_1 ~\gamma~ \ldots ~\gamma~ x_h$ from the left \\

(2.11)~~~ $ \begin{array}{l}
               x_1 ~\gamma~ x_2 ~\approx_{\cal A}~ \alpha ( x_1, x_ 2 ) ~\approx_{\cal A}~ a x_1 + b x_2 \\
               x_1 ~\gamma~ x_2 ~\gamma~ x_3 ~\approx_{\cal A}~ ( a x_1 + b x_2 ) ~\gamma~ x_3
               ~\approx_{\cal A}~ a^2 x_1 + a b x_2 + b x_3 \\
               x_1 ~\gamma~ x_2 ~\gamma~ x_3 ~\gamma~ x_4 ~\approx_{\cal A}~ ( a^2 x_1 + a b x_2 + b x_3 ) ~\gamma~ x_4 ~\approx_{\cal A}~ \\
               ~~~~~~~~~~~~\approx_{\cal A}~ a^3 x_1 + a^2 b x_2  + a b x_3 + b x_4 \\
               \vdots \\
               x_1 ~\gamma~ x_2 ~\gamma~ \ldots ~\gamma~ x_h ~\approx_{\cal A}~
               a^{h-1} x_1 + a^{h-2} b x_2 + \ldots + a b x_{h-1} + b x_h \\
               \vdots
            \end{array} $ \\

and alternatively, by successively grouping the terms in $x_1 ~\gamma~ \ldots ~\gamma~ x_h$ from the right, we obtain \\

(2.12)~~~ $ x_1 ~\gamma~ x_2 ~\gamma~ \ldots ~\gamma~ x_h
~\approx_{\cal A}~
               a x_1 + a b x_2 + \ldots + a b^{\,h-2} x_{h-1} + b^{\,h-1} x_h $ \\

Furthermore, we can apply (1.21) to (2.11), (2.12), and obtain a lot of other equivalences. \\

For convenience, however, let us assume $a = b \geq 2$, then (2.11), (2.12) give \\

(2.13)~~~ $ x_1 ~\gamma~ x_2 ~\gamma~ \ldots ~\gamma~ x_h
~\approx_{\cal A}~
               a x_{i_1} + a^2 x_{i_2} + \ldots + a^{h-1} x_{i_{h-1}} + a^h x_{i_h} $ \\

where $i_1, \ldots , i_h$ is any permutation of $1, \ldots , h$. \\

Now, since $\alpha$ in (2.9) is {\it not} associative, it follows in
view of Theorem 2.1., that the mapping $\iota$ in (2.5) is
surjective, but {\it not}
injective. \\

Let us therefore see how much $\bigotimes_{\cal A} X$ is {\it smaller} than $X$. \\

This means that, for given $x, y \in X,~ x \neq y$, we have to find out when nevertheless $\iota ( x ) = \iota ( y )$, which means \\

(2.14)~~~ $ x ~\approx_{\cal A}~ y $ \\

In view of (2.10), (2.13), we have therefore \\

(2.15)~~~ $ \begin{array}{l}
                    x = a x_1 + \ldots + a^h x_h ~\approx_{\cal A}~ x_1 ~\gamma~ \ldots ~\gamma~ x_h ~\approx_{\cal A}~ \\
                    ~~~~~\approx_{\cal A}~ y_1 ~\gamma~ \ldots ~\gamma~ y_k ~\approx_{\cal A}~ a y_1 + \ldots + a^k y_k = y
             \end{array} $ \\

for certain $h, k \in \mathbb{N},~ h, k \geq 2$ and $x_1, \ldots , x_h, y_1, \ldots , y_k \in X$. \\

On the other hand, it is obvious that \\

(2.16)~~~ $ 0 < x = a x_1 + \ldots + a^h x_h ~~~\Longrightarrow~~~ x \geq a \geq 2 $ \\

Thus (2.14) - (2.16) give \\

(2.17)~~~ $ x, y \in X,~ x, y < a,~ x ~\approx_{\cal A}~ y ~~~\Longrightarrow~~~ x = y $ \\

In other words, the mapping, see (2.5) \\

(2.18)~~~ $ X \ni x \stackrel{\iota}\longmapsto \bigotimes_{\cal A} x \in \bigotimes_{\cal A} X $ \\

is {\it injective} at least on $\{ 0, 1, 2, \ldots , a - 1 \} \subset \mathbb{N} = X$. \\ \\

{\bf 3. Beyond Classical Superpositions} \\

A special interest in the extended concept of superposition
presented above is in the fact that, as seen in the sequel, it
incorporates as a particular
case the non-classical quantum type superposition. \\

First however, for the sake of clarity, we consider superposition in a somewhat more general situation, that that in the usual quantum case. \\

{\bf Non-Classical Superpositions} \\

Let $X$ and $K$ be two arbitrary sets. We shall consider
superposition on the set $Y = K \times X$, and which is defined by a
family ${\cal A}$ of
$(n+m)$-ary relations on $Y$ that consists this time of one single $(2+1)$-ary relation $A \subseteq Y^2 \times Y$. In this case (1.22) becomes \\

(3.1)~~~ replace $( c_1, x_1 ) ~\gamma~ ( c_2, x_2 ) ~\gamma ~ ( d,
y ) ~\gamma~ \ldots $

           ~~\hspace*{1cm} with $( c, x ) ~\gamma~ ( d, y ) ~\gamma~ \ldots $

           ~~\hspace*{1cm} or vice-versa, where

           ~~\hspace*{1cm} $( ( c_1, x_1 ), ( c_2, x_2 ), ( c,  y ) ) \in A $ \\

{\bf Quantum Superpositions} \\

An example of the above situation comes from {\it quantum
superposition} in which case $X$ is a given Hilbert space $H$, while
$K$ is the field $\mathbb{C}$ of complex numbers, thus the elements
$( c, x ) \in Y = K \times X = \mathbb{C} \times H$ are of the form
$( c, | \psi > )$, where $c \in
\mathbb{C}$ and $|\psi >\, \in H$. \\
In this quantum case one excludes elements $( c, x ) = ( c, | \psi >
)$, with $c = 0 \in \mathbb{C}$ or $| \psi >\, = 0 \in H$, and
furthermore, one identifies any two elements $( c, x ) = ( c, | \psi
> )$ and $( c\,', x ) = ( c\,', | \psi > )$, as long as $c, c\,'
\neq 0$. As it turns out, however, the above construction involved
in superposition can be done without these two restrictions, which
can of course be brought into consideration at the end
of the construction. \\

In this quantum case, we can specify the way the $(2+1)$-ary relation $A \subseteq Y^2 \times Y$ is chosen, namely \\

(3.2)~~~ $ A = \left \{~ ( ( c_1, x_1 ), ( c_2, x_2 ), ( c,  y ) )
\in Y^3 ~
                                     \begin{array}{|l}
                                        ~1)~ c_1, c_2, c \neq 0 \\ \\
                                        ~2)~ x_1, x_2 \neq 0 \\ \\
                                        ~3)~ y = d_1 x_1 + d_2 x_2, ~\mbox{where}~ d_1, d_2 \in K \\ \\
                                        ~~~~~~\mbox{and}~ | d_1 |^2 || x_1 ||^2 + | d_2 |^2 || x_2 ||^2 = 1
                                     \end{array} ~\right \} $ \\

thus (3.1) takes the particular form \\

(3.3)~~~ replace $( c_1, x_1 ) ~\gamma~ ( c_2, x_2 ) ~\gamma ~ ( d,
y ) ~\gamma~ \ldots $

           ~~\hspace*{1cm} with $( c, x ) ~\gamma~ ( d, y ) ~\gamma~ \ldots $

           ~~\hspace*{1cm} or vice-versa, where \\

           ~~\hspace*{1cm} 1)~ $ c_1, c_2, c \neq 0 $ \\

           ~~\hspace*{1cm} 2)~ $ x_1, x_2, x \neq 0 $ \\

           ~~\hspace*{1cm} 3)~ $ x = d_1 x_1 + d_2 x_2 $ \\

           ~~\hspace*{1cm} 4)~ $ | d_1 |^2 || x_1 ||^2 + | d_2 |^2 || x_2 ||^2 = 1 $ \\ \\

{\bf Appendix : Free Semigroups, etc.} \\ \\

We recall for convenience several basic concepts and constructions related to semigroups, [3,7]. \\ \\

{\bf A1. Semigroups} \\

A {\it semigroup} is a structure $( S, \star )$, where $S$ is a
nonvoid set and $\star : S \times S \longrightarrow S$ is a binary
operation on $S$ which
is {\it associative}, that is, it satisfies the condition \\

(A1.1)~~~ $ \star ( u, \star ( v, w ) ) = \star ( \star ( u, v ), w ),~~~ u, v, w \in S $ \\

Here we recall that it is customary to denote $\star ( u, v )$
simply by $u \star v$, for $u, v \in S$. Consequently, the above
associativity condition
is equivalent with \\

(A1.2)~~~ $ u \star ( v \star w ) = ( u \star v ) \star w,~~~ u, v, w \in S $ \\

We note that semigroups need not always have to be commutative, or have neutral elements. \\

Given any semigroup $( S, \star )$, we associate with it the semigroup $( S^1, \star )$, where \\

(A1.3)~~~ $ S^1 = \begin{array}{|l}
                                  S ~~\mbox{if}~ S ~\mbox{has a neutral element}~ e \\ \\
                                  S \cup \{ e \} ~~\mbox{if}~ S ~\mbox{does not have a neutral element}
                     \end{array} $ \\

In the second case, we extend $\star$ to $S^1 \times S^1$ in the
obvious manner, namely, by defining $u \star e = e \star u = u$, for
$u \in S$. And in
this case $( S^1, \star )$ will be a semigroup with the neutral element $e$. \\

A fundamental concept which can relate semigroups to one another is
introduced now. Given two semigroups $( S, \star )$ and $( T,
\diamond )$, a mapping $f :
S \longrightarrow T$ is called a {\it homomorphism}, if and only if \\

(A1.4)~~~ $ f ( u \star v ) = f ( u ) \diamond f ( v ),~~~ u, v \in S $ \\

If such a homomorphism $f$ is injective, then it is called a {\it
monomorphism}. In case a homomorphism $f$ is surjective, then it is
called an {\it epimorphism}. And if a homomorphism $f$ is surjective
and it also has an inverse mapping $f^{-1} : T \longrightarrow S$
which is again a homomorphism, then it is called an {\it
isomorphism}. It follows that isomorphic semigroups are identical
for all purposes, and they only differ in the notation of
their elements, or of the semigroups themselves. \\ \\

{\bf A2. Free Semigroups} \\

Free semigroups are of fundamental importance since, as indicated in
Proposition A2.1. below, all semigroups can be obtained from them in
a natural
manner. \\

Given any nonvoid set $E$, we denote by $E^+$ the set of all finite sequences \\

(A2.1)~~~ $ a_1, a_2, \ldots , a_n $ \\

where $n \geq 1$ and $a_i \in E$, with $1 \leq i \leq n$. Further, we define on $E^+$ the binary operation $\ast$ as follows \\

(A2.2)~~~ $ ( a_1, a_2, \ldots , a_n ) \ast ( b_1, b_2, \ldots , b_m
) =
                                  a_1, a_2, \ldots , a_n, b_1, b_2, \ldots , b_m $ \\

that is, simply by the juxtaposition or concatenation of sequences in  $E^+$. Then it is easy to see that \\

(A2.3)~~~ $ ( E^+, \ast ) $ is a semigroup \\

since $\ast$ is obviously associative. However, if $E$ has at least
two elements, then $\ast$ is clearly not commutative. Also, the
semigroup
$( E^+, \ast )$ does not have a neutral element, regardless of the number of elements in $E$. \\

A consequence of the associativity of $\ast$ is that the elements of
the semigroup $E^+$ can be written in the following form which is
alternative to
(A2.1), namely \\

(A2.4)~~~ $ a_1 \ast a_2 \ast \ldots \ast a_n $ \\

where $n \geq 1$ and $a_i \in E$, with $1 \leq i \leq n$. \\

The semigroup $( E^+, \ast )$ is called the {\it free semigroup on}
$E$, and the meaning of that term will result from the two important
properties
presented next. \\

First, the mapping \\

(A2.5)~~~ $ E \ni a \stackrel{i}\longmapsto a \in E^+ $ \\

is {\it injective}, and obviously, it is never surjective. \\

Second, the free semigroup $( E^+, \ast )$ on $E$ has the following
{\it universality property}. Given any semigroup $( S, \star )$ and
any mapping
$j : E \longrightarrow S$, there exists a {\it unique} homomorphism $f : E^+ \longrightarrow S$ such that \\

(A2.6)~~~ $ j = f \circ i $ \\

or in other words, the diagram commutes

\begin{math}
\setlength{\unitlength}{1cm} \thicklines
\begin{picture}(13,6)

\put(0,2.4){$(A2.7)$} \put(1.1,5){$E$} \put(2,5.1){\vector(1,0){7}}
\put(9.5,5){$E^+$} \put(5.5,5.4){$i$}
\put(1.7,4.5){\vector(1,-1){3.3}} \put(9.5,4.5){\vector(-1,-1){3.3}}
\put(5.4,0.6){$S$} \put(2.8,2.5){$j$} \put(8.1,2.5){$\exists !~~f$}

\end{picture}
\end{math} \\

The above universality property of semigroups has an important immediate consequence, namely \\

{\bf Proposition A2.1.} \\

Every semigroup is the homomorphic image of a free semigroup. \\

{\bf Proof.} \\

Let $( S, \star )$ be a semigroup, then we can take $E = S$ and $j =
id_S$ in (A2.7), and obtain the commutative diagram

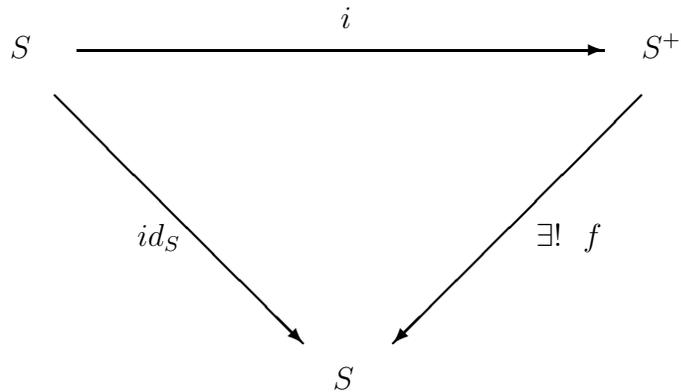
\begin{math}
\setlength{\unitlength}{1cm} \thicklines
\begin{picture}(13,7)

\put(1.1,5){$S$} \put(2,5.1){\vector(1,0){7}} \put(9.5,5){$S^+$}
\put(5.5,5.4){$i$} \put(1.7,4.5){\vector(1,-1){3.3}}
\put(9.5,4.5){\vector(-1,-1){3.3}} \put(5.4,0.6){$S$}
\put(2.8,2.5){$id_S$} \put(8.1,2.5){$\exists !~~f$}

\end{picture}
\end{math} \\

However, in view of (A2.6), obviously $f$ is surjective. Thus $S$ is indeed the homomorphic image of the free semigroup $S^+$. \\ \\

{\bf A3. Quotient Constructions} \\

Let $E$ be a nonvoid set and $\approx$ an equivalence relation on $E$. Then the {\it quotient set} \\

(A3.1)~~~ $ E / \approx $ \\

is defined as having the elements given by the {\it cosets} \\

(A3.2)~~~ $ ( a )_\approx = \{~ b \in E ~~|~~ b \approx a ~\},~~~ a \in E $ \\

thus each coset $( a )_\approx$ is the set of all elements $b \in E$ which are equivalent with $a$ with respect to $\approx$. The coset $( a )_\approx$ is also called the equivalence class of $a$ with respect to the equivalence relation $\approx$. It follows that the mapping \\

(A3.3)~~~ $ i_\approx : E \ni a \longmapsto ( a )_\approx \in E / \approx $ \\

is surjective, and it is called the {\it canonical quotient mapping}. \\

A useful way to obtain equivalence relations on any given set $E$ is
through the construction called {\it transitive closure}. Namely,
given any family $( \equiv_i )_{i \in I}$ of symmetric binary
relations on $E$, then we define the equivalence relation $\approx$
on $E$ as follows. If $a, b \in E$,
then \\

(A3.4)~~~ $ a \approx b $ \\

holds, if and only if $a = b$, or there exist $c_0, c_1, c_2, \ldots , c_n \in E,~ i_1, i_2, i_3, \ldots , i_n \in I$, such that \\

(A3.5)~~~ $ a ~=~ c_0 ~\equiv_{i_1}~ c_1 ~\equiv_{i_2}~ c_2 ~\equiv_{i_3}~ \ldots ~\equiv_{i_n}~ c_n ~=~ b $ \\

An alternative and equivalent way to construct quotient spaces is
through {\it partitions}. Given a partition of $E$ by the family
${\cal E} = ( E_i )_{i \in I}$ of subsets of $E$. Then one can
associate with it an equivalence relation
$\approx_{\cal E}$ on $E$, defined for $a, b \in E$, by \\

(A3.6)~~~ $ a \approx_{\cal E} b ~~~\Longleftrightarrow~~~ \exists~~ i \in I ~:~ a, b \in E_i $ \\

Obviously, in this case we have for $a \in E$ and $i \in I$ \\

(A3.7)~~~ $ a \in E_i ~~~\Longleftrightarrow~~~ E_i = ( a )_{\approx_{\cal E}} $ \\

in other words, the equivalence class $( a )_{\approx_{\cal E}}$ of $a$ with respect to $\approx_{\cal E}$ is precisely the set $E_i$ in the partition ${\cal E}$ to which $a$ belongs. Consequently \\

(A3.8)~~~ $ E / \approx_{\cal E} ~=~ \{~ E_i ~~|~~ i \in I ~\} $ \\ \\

{\bf A4. Congruences} \\

Let $( S, \star )$ be any semigroup. An equivalence relation
$\approx$ on $S$ is called a {\it congruence} on $( S, \star )$, if
and only if it is
compatible with the semigroup operation $\star$ in the following sense \\

(A4.1)~~~ $ u \approx v ~~~\Longrightarrow~~ u \star w \approx v \star w,~~ w \star u \approx w \star v $ \\

for all $u, v, w \in S$. \\

The importance of such a congruence is that the resulting quotient $S / \approx $ of $S$ leads again to a {\it semigroup}, namely \\

(A4.2)~~~ $ ( S, \star ) / \approx ~~~=~~~ (\, S / \approx,\, \diamond ) $ \\

where the binary operation $\diamond$ on $S / \approx$ is defined by \\

(A4.3)~~~ $ ( u )_\approx \diamond ( v )_\approx = ( u \star v )_\approx,~~~ u, v \in S $ \\

also, the canonical quotient mapping, see (A3.4) \\

(A4.4)~~~ $ S \ni u \longmapsto ( u )_\approx \in S / \approx $ \\

is a {\it surjective homomorphism}, thus an {\it epimorphism}. \\

Furthermore, let $( S, \star )$ and $( T, \diamond )$ be two
semigroups and $f : S \longrightarrow T$ a morphism between them.
Then the binary relation
on $S$ given by \\

(A4.5)~~~ $ ker f = \{~ ( u, v ) \in S \times S ~~|~~ f ( u ) = f ( v ) ~\} $ \\

is a congruence on $( S, \star )$, and there exists a monomorphism $g : ( S, \star ) / ker f \longrightarrow ( T, \diamond )$, such that, see (A3.3) \\

(A4.6)~~~ $ f = g \circ i_{ker f} $ \\

which means that the diagram commutes \\

\begin{math}
\setlength{\unitlength}{1cm} \thicklines
\begin{picture}(13,6)

\put(0,2.5){$(A4.7)$} \put(1.1,5){$S$}
\put(2,5.1){\vector(1,0){7.3}} \put(9.9,5){$T$} \put(5.5,5.4){$f$}
\put(1.7,4.5){\vector(1,-1){3.3}} \put(6.5,1.3){\vector(1,1){3.3}}
\put(5.1,0.6){$S / ker f$} \put(2.5,2.5){$i_{ker f}$}
\put(8.1,2.5){$g$}

\end{picture}
\end{math} \\

Given now a partition ${\cal S} = ( S_i )_{i \in I}$ of $S$, then in view of (A3.6), it leads to an equivalence relation $\approx_{\cal S}$ on $S$. Now in view of (A4.1), (A3.6), (A3.7), it is obvious that the equivalence relation $\approx_{\cal S}$ on $S$ will be a {\it congruence} on $( S, \star )$, if and only if, for every $i \in I,~ u, v \in S_i,~ w \in S$, we have \\

(A4.8)~~~ $ ( u \star w )_{\approx_{\cal S}} = ( v \star w
)_{\approx_{\cal S}},~~~
               ( w \star u )_{\approx_{\cal S}} = ( w \star v )_{\approx_{\cal S}} $ \\ \\

{\bf A5. On the Structure of Semigroups} \\

A basic structural result on semigroups was presented in Proposition
A2.1. Here we mention another one which brings a new and additional
such
structural insight. \\

Given any nonvoid set $E$, we denote by ${\cal T}_E$ the set of all
mappings $f : E \longrightarrow E$ of $E$ into itself. It is easy to
see that $( {\cal T}_E, \circ )$ is a semigroup, where $\circ$ is
the usual composition of mappings, and it is called the {\it full
transformation semigroup} on
$E$. \\

A classic extension of Cayley's theorem for groups is presented in, see (A1.3) : \\

{\bf Proposition A5.1.} \\

Given any semigroup $( S, \star )$. Then there is a monomorphism,
that is, an injective morphism $f : S \longrightarrow {\cal
T}_{S^1}$.

\hfill $\Box$ \\

In other words, every semigroup is the subsemigroup of a full
transformation semigroup. And in case the semigroup $( S, \star )$
has a neutral element,
then in view of (A1.3), it is simply the full transformation semigroup of $S$ into itself. \\ \\

{\bf A6. A Commutative Semigroup} \\

Let $( E^+, \ast )$ be the free semigroup generated by the nonvoid set $E$, see (A2.3). Given two sequences, see (A2.4) \\

(A6.1)~~~ $ a_1 \ast \ldots \ast a_n,~~~ b_1 \ast \ldots \ast b_m \in E^+ $ \\

we define \\

(A6.2)~~~ $ a_1 \ast \ldots \ast a_n ~\approx~ b_1 \ast \ldots \ast b_m $ \\

if and only if the two sequences are the same, or differ by a
permutation of their elements. Then obviously $\approx$ is an
equivalence relation on $E^+$
which is also a congruence on $( E^+, \ast )$. It follows that \\

(A6.3)~~~ $ ( E^+, \ast ) / \approx $ \\

is a commutative semigroup. \\

\end{document}